\documentclass{article}[12pt]
\usepackage{color,times,amsmath,amsfonts,latexsym,epsfig,epsf,colordvi}
\usepackage{url,hyperref}
\usepackage[english]{babel}

\newcommand{\CC}{\mathbb {C}}
\newcommand{\RR}{\mathbb {R}}

\newcommand{\bp}{\begin{pmat}}
\newcommand{\ep}{\end{pmat}}

\author{  Frank Uhlig \thanks{Department of Mathematics and Statistics, Auburn 
University, Auburn, AL 36849-5310 \ (uhligfd@auburn.edu)}}

 \title{\vspace*{-8mm} Coalescing Eigenvalues and Crossing Eigencurves of 1-Parameter Matrix Flows\\[0mm]}

\begin{document}
\date{~}
\thispagestyle{empty}
\maketitle

\thispagestyle{empty}

\vspace*{-14mm}
\begin{center} { \bf Abstract  } \\[2mm]
\begin{minipage}{150mm}
We investigate the eigenvalue curves of 1-parameter hermitean and general complex or real matrix flows $A(t)$ in light of their geometry and the uniform decomposability of $A(t)$ for all parameters $t$. The often misquoted and misapplied results by Hund and von Neumann and by Wigner for eigencurve crossings from the late 1920s are clarified for hermitean matrix flows $A(t) = (A(t))^*$. A  conjecture on extending these results to general non-normal or non-hermitean 1-parameter matrix flows is formulated and investigated. An algorithm to compute the block dimensions of  uniformly decomposable hermitean matrix flows is described and tested. The algorithm uses the ZNN method to compute the time-varying matrix eigenvalue curves of $A(t)$ for $t_o \leq t\leq t_f$. Similar efforts for general complex matrix flows are described. This extension  leads to many new and open problems. Specifically, we point to the difficult relationship between the geometry of eigencurves for general complex matrix flows $A(t)$ and a general flow's decomposability into blockdiagonal form via one fixed unitary or general matrix similarity for all parameters $t$.
\end{minipage}\\[-1mm]
\end{center}  
\thispagestyle{empty}

\noindent{\bf Keywords:}  matrix  eigenvalues, time-varying matrix flows, eigenvalue curve, eigencurve crossing, Hund-von Neumann-Wigner Theorem, Zhang Neural Network,    numerical matrix algorithm, decomposable matrix, block diagonal matrix\\[-4mm]

\noindent{\bf AMS :} 15A60, 65F15, 65F30, 15A18\\[-7mm]

\pagestyle{myheadings}
\thispagestyle{plain}
\markboth{Frank Uhlig}{Eigencurves of  1-parameter matrix flows }


\section{Introduction }\vspace*{-2mm} 

The eigenvalues of real or complex matrix flows $A(t)_{n,n}$ have been studied for more than 90 years, see \cite{FH1927}, \cite{NW} for the earliest papers that came about within foundational quantum theory. Parameter-dependent eigenvalue curves of hermitean matrix flows and their possible crossings  have become important for studies on stability and bifurcation such as for molecular aspects of quantum and chemical physics, in the study of energy surfaces, in structural analyses, in antenna theory and  in other areas, see e.g. \cite{KST81} or \cite{Sch16} and specifically  \cite[p. 520]{DPP13} for a  listing of recent references. In the early 1900s it was most important for quantum mechanics to understand  whether parameter varying matrix eigenvalue curves would  intersect or cross, leading to different practical results and quantum state implications if they would. The two fundamental papers \cite{FH1927}, \cite{NW} 'proved' that hermitean matrix flows that depend on one   (or two)  parameters, such as $A(t) = F + t G$ or $A(t) = F + t G+t^2H$  for constant hermitean matrices $F, G$, and/or $H$ would not allow eigencurve crossings. \\[1mm]
This classical 'result' has not been challenged and has unfortunately been repeated (see e.g. \cite[p. 519, 2nd paragraph in Introduction]{DPP13}) dozens of times in the literature until 2018 when normal non-hermitean matrices were named as exceptions to the von Hund-Neuman-Wigner no-crossing rule in \cite[Ex. 7.1, Fig. 7.2, p. 1739, 1740]{LM} in computational studies of the field of values of a constant matrix. As the simplest and worst case counter example for the classical result with hermitean 1-parameter matrix flows, note that if a $2n$ by $2n$ 1-parameter matrix flow $A(t) = F(t) + t G(t)$ is generated from two compatibly dimensioned  block diagonal hermitean matrix flows $F(t) = diag(f,f)_{2n,2n}$ and $G(t) = diag(g,g)_{2n,2n}$ with  repeated hermitean $n$ by $n$ blocks $f(t)_{n,n}$ and $g(t)_{n,n}$, respectively, then every point on every  eigencurves of $A(t)_{2n,2n}$ is doubly covered, i.e., every point of every eigencurve of $A(t)_{2n,2n}$ is a 'crossing point'. \\[1mm] 
One  explanation for this historical oversight by the 'eigenvalue crossing' community is the fact that a simple blockdiagonal counter example would never appear in  or become an issue in quantum physics nor in the study of single atoms or molecules.\\
 But for mathematicians this simple counterexample produces many note-worthy challenges: can there be  eigenvalue curve crossings in more general 1-parameter matrix flow settings?  When do they occur, if ever? How can they be found? Which is the coarsest block-diagonalisation for a given 1-parameter general complex or  hermitean matrix flow $A(t)$, which the finest? How can the decomposition block sizes for a decomposable matrix flow $A(t)$ be determined from its eigencurves? Can an actual  block diagonalization of a decomposable matrix flow $A(t)$ be computed?\\[1mm]
The initial paper  \cite{NW} by von Neumann and Wigner studied 1- to 3-variable dependent hermitean matrix flows in regards to eigencurve crossings. Recent work of Dieci et al. \cite[Ex. 4.1, 4.2, p. 533 - 535]{DPP13} has created an algorithm for computing multi-variable  eigencrossing data points  for 1- to 3-parameter hermitean matrix flows  of the form $f(x,y,z)F + g(x,y,z)G + h(x,y,z)H$ with constant matrices $F,G,H$ and multi-nomials $f, g$, and $h$.\\[2mm]
Here we consider 1-parameter hermitean and general complex or real matrix flows $A(t) = (a_{k,j}(t))_{n,n}$ with 1-parameter varying entry functions. This extends the earlier $A(t) = F + tG \ (+\ t^2H)$ approach in two directions, namely 
to  general matrix flows with  general entry functions.\\[1mm]
Section 2 deals with 1-parameter hermitean matrix flows $A(t)$ and eigencurve crossings and relates the latter to equivalent separable or decomposable block matrix flows, where the term 'equivalent' means uniformly  via one  fixed matrix  similarity $S^{-1}A(t)S$ for all parameters $t$. An algorithm for finding the block dimensions of decomposable hermitean flows is given and open questions are raised. Section 3 deals with general complex or real 1-parameter matrix flows. Throughout we use the Zhang Neural Network (ZNN) method for the eigenanalyses of time-varying  matrix flows as developed in \cite{ZYLUH,YZMeZNN} .\\[-8mm]

\section{Hermitean Matrix Flows and  Eigencurve Crossings}

This section deals with 1-parameter matrix flows $A(t)_{n,n}$, their eigenvalue curves for a time or parameter interval $t_o \leq t \leq t_f$,  and the notion of matrix separability or matrix decomposition. More specifically we deal  with hermitean matrix flows in this section; general matrix flows and their eigencurves are discussed in Section 3.\\[-3mm]

\noindent
{\bf Definition  :}\\[1mm]
\hspace*{8mm} \begin{minipage}{150mm}  {\bf (1)} A constant square matrix $A_{n,n}$ is called \emph{separable} or \emph{decomposable} if $A$ is similar to a proper blockdiagonal matrix,. Here and below 'proper' means that $A$'s blockdiagonal representation has at least two diagonal blocks.\\
{\bf (2)} A 1-parameter real or complex square matrix flow $A(t)_{n,n}$ is called \emph{structurally separable} or \emph{structurally decomposable} on an interval $[t_o, t_f] \subset \RR$ if each $A(t)$ can be reduced uniformly to the same  proper blockdiagonal form via the same fixed matrix similarity. 
\end{minipage}

\vspace*{1mm}

\noindent
Note that an indecomposable or  decomposable matrix flow might contain specific matrices $A(t)$ that may be  reduced further for a specific value of $t \in [t_o, t_f]$  if, for example, some strategic entries in the common block diagonal form of all $A(t)$ become zero at some $t$. \\[1mm]
Obviously  the eigencurves of a block diagonal hermitean matrix flow $A(t) = diag(A_1, A_2, ..., A_k)$  are simply the superpositions of the eigencurves of each of its individual  matrix blocks $A_i$. If $k = 2$ and the eigencurves of the first diagonal block $A_1$ hover around 100 on a given interval $[t_o,t_f]$ and those of the second  block $A_2$ hover around --50 in value, for example, then there will likely be no crossings among the set of eigencurves for $A(t)$ on the given interval and then the eigencurves  carry little information regarding  the possible decomposability of $A(t)$.  In this case it might be advisable to enlarge the interval $[t_o, t_f]$ since the decomposability of parameter-varying matrix flows is a global property. We will study and learn more about this phenomenon later on.
 If on the other hand there are observed eigencurve crossings for a hermitean matrix flow $A(t)_{n,n}$, then $A$ is separable or decomposable, i.e., there is a nonsingular constant matrix $S_{n,n}$ so that $S^{-1}\cdot A(t) \cdot S$ is uniformly and properly block diagonal for all $t$ as we shall see. Throughout we will use unitary similarities  $ S = U \text{ with } UU^* = I $  in order to not affect the eigenvalue conditioning of the flows. All our programs are designed to work with ordinary similarities and with hermitean or orthogonal ones. \\[1mm]
 Note that once we know the eigencurve data for $A(t)_{n,n}$ approximately we can interpolate the eigenvalue curves  and form a diagonal matrix flow $\tilde A(t) = diag(a_1(t), ..., a_n(t))$ with the individual eigencurve functions $a_i(t)$ of $A(t)$ in successive diagonal positions. Therefore every matrix flow might potentially come from a completely decomposable, i.e., a diagonal matrix flow. Therefore looking for the finest possible decomposition structure of matrix flows is futile. It would be   more sensibly to try and find the coarsest proper decomposition structure of  a given  1-parameter varying matrix flow instead.  This is a new and worthwhile question that combines function geometry with matrix analysis.\\[1mm]
The classical hermitean flow results of Hund \cite{FH1927} and von Neumann and Wigner \cite{NW}, stated correctly for indecomposable matrix flows, are as follows.\\[1mm]
{\bf Hund-von Neumann-Wigner Theorem \ \cite{FH1927,NW} :} (abbreviated by HvNW) \\
If $A(t)$ is an indecomposable 1-parameter hermitean matrix flow, then\\[1mm] \enlargethispage{20mm}
\hspace*{3mm} {\bf (a}) the eigenvalue curves in $\RR^2 \simeq \CC$ of $A(t)$ do not intersect \ and \\[1mm]
\hspace*{3mm} {\bf (b}) if two eigenvalue curves approach each other, they avoid crossing each other by veering off in a hyperbolic \linebreak
\hspace*{9.2mm}way where the approaching angle of either eigencurve equals the leaving angle of the other eigenvalue curve \linebreak
\hspace*{9.2mm}after their close encounter.\\[1mm]
For a proof of part b) using Schr\"odinger's perturbation method and asymptotics, see \cite[Section 2, p. 469]{NW} and its Figure 1 copied below. A similar process of veering off was   observed and established for multi-parameter eigencurves in \cite{DPP13} and termed 'conical intersection of eigenvalues' there.\\[-6mm]
\begin{center}
\includegraphics[width=51mm]{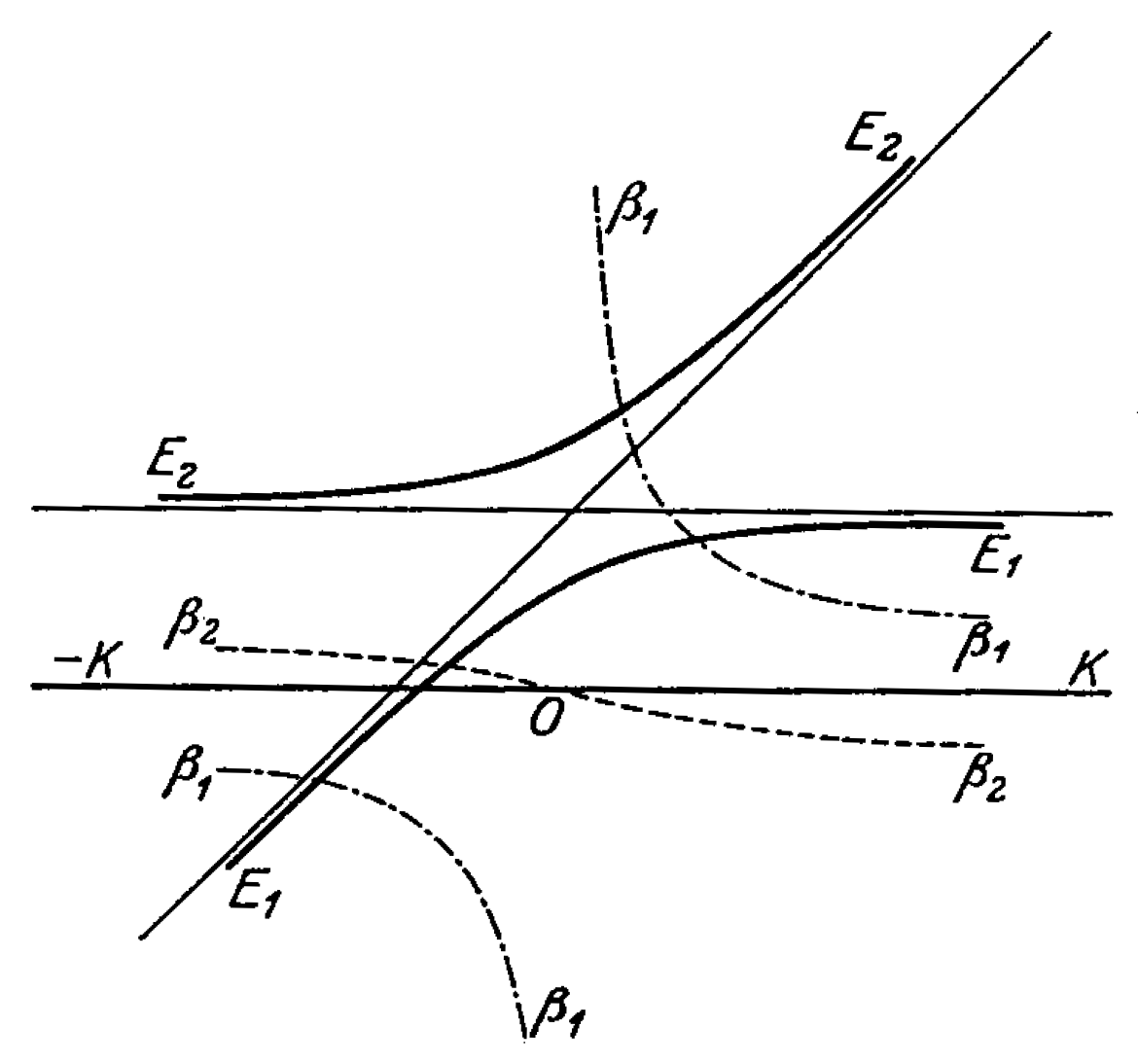} \\[0mm]
Figure 1 \ \ \ \ \ (from  \cite[p. 469]{NW}) 
\end{center} 
\vspace*{-1mm}
Throughout we  use the fast and highly accurate Zhang Neural Network approach, see \cite{ZYLUH,YZMeZNN,FUZhFoV},
 to plot the eigencurves of 1-parameter varying matrix flows $A(t)$ and determine  increasingly finer coarse block-diagonal decomposition sizes that can be obtained for a given matrix flow.\\[1mm]
In Figure 2 we plot the real eigencurves  for a hermitean time-varying  random entry matrix flow $A(t)_{11,11}$.\\[-6mm] 
\begin{center}
\includegraphics[width=110mm]{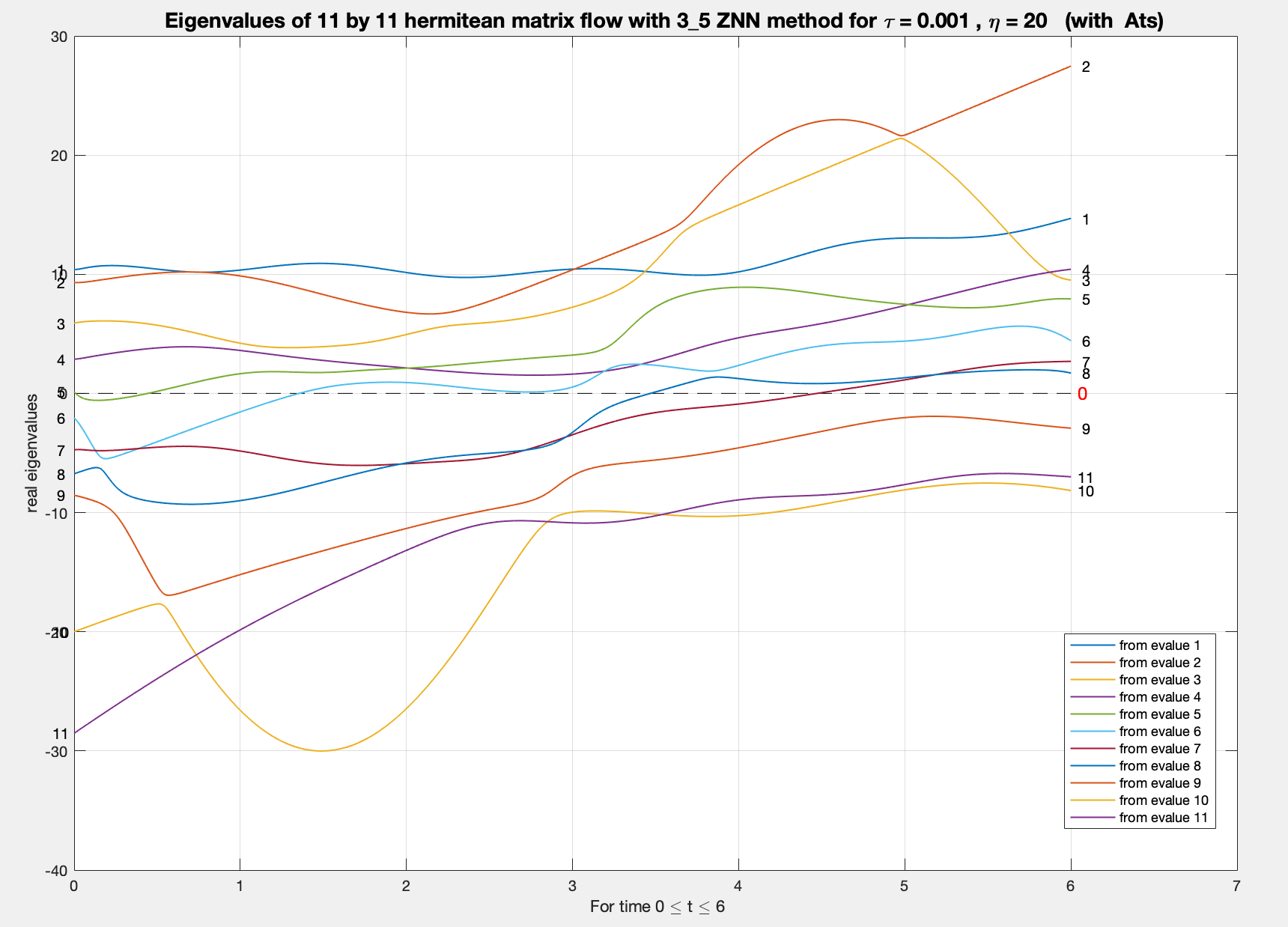} \\[0mm]
Figure  2
\end{center}
\newpage

\noindent
In Figure 2, the 11 real eigencurves for $A(t)$ are indexed by their eigenvalues at $t_o = 0$ in descending order along the left edge. They are traced in differing colors until $t_f=6$ where their curve numbers are repeated for clarity. A colored legend panel  -- if color is available -- further helps with determining which eigencurve crosses which. We compute eigencurve crossings in our Matlab code {\tt Chermitmatrixfloweig.m} \footnote{Matlab codes for all main and auxiliary m-files are collected and referenced in \cite{FUEigencurves}.}  which stores the crossings data in  an $n-1$ by $n+1$  matrix R1 for $n$ by $n$ hermitean matrix flows. For   Figure 2, R1's first 5 columns are \\[1mm]
\hspace*{56mm}\begin{tabular}{c|cccc}
curve & \multicolumn{4}{c}{crosses eigencurve}\\
number& \multicolumn{4}{c}{ with label}\\ \hline
1&2&3&0& 0\\
2&0&0&0&0\\
3&4&0&0&0\\
4&5&6&0&0\\
5&0&0&0&0\\
6&7&0&0&0\\
7&8&0&0&0\\
8&0&0&0&0\\
9&0&0&0&0\\
10&11&0&0&0
\end{tabular} \ .\\[2mm]
The code computes  the vector $ve =  (1, -1, -1, 1, -1, -1, 1, -1, 2, 3, -3) \in \RR^{11}$ from  R1. This integer vector $ve$ separates the eigencurves into crossing groups which fall into 5 sets here as indicated by $ve$'s five distinct entries 1, --1, 2, 3, and --3 here. The vector $ve$ always starts with a 1 for the first, the topmost starting curve and all  other entries are  initially set to zero. All eigencurves  that  the first eigencurve crosses are labeled with --1 in $ve$. If there are   zero entries in $ve$ afterwards,  the algorithm starts from there for  the indices  2  and --2 and indexes  the  eigencurves below in R1 as before. This process is repeated until all 11 eigencurves are labelled  with positive or negative integers in succession in $ve$. If there are data clashes where a nonzero entry $k$ in $ve$ does not conform with a new crossing requirement of $-k$ in a position, this indicates the need for an additional eigencurve group. Therefore the integer label for this entry of $ve$ is upped to $k+1$, all subsequent $ve$ entries in this row of R1 are reset to zero. And the remaining crossing numbers  for the current curve in R1 are skipped for a fresh start, and we continue with the next row of R1. Eigencurves with opposite signed integer values in $ve$ cross each other. Further details are listed in the \%  code line comments of  {\tt  Chermitmatrixfloweig.m}.\\
We advise to occasionally try and create the decomposition data vector $ve$ for a given hermitean matrix flow by hand, using pencil and paper, the  crossing matrix R1 and  curve plot verifications to learn how the automatic code uses part (a) of the HvHW Theorem.\\
 In our chosen example and according to the currently available $ve=  (1, -1, -1, 1, -1, -1, 1, -1, 2, 3, -3)$, there are two  diagonal blocks of $A(t)_{11,11}$ associated with the number 1: three eigencurves carry a 1 and five a --1, meaning that the ones with the label 1 cross the ones with the label --1 and vice versa, while no others are crosses  by  or crossing these 8 eigencurves. Moreover there is a single, the 9th eigencurve associated with the label 2 in $ve$, that is not crossing those above, nor those below as can readily be seen in R1 and in Figure 2. It is separate from, i.e., non-crossing  all other eigencurves. Finally there are two eigencurves that cross  one another and they were given the labels 3 and --3 in positions 10 and 11 of $ve$.\\
  Thus far we have processed the eigencurve crossing data in R1 according only to part (a) of the HvNW Theorem.\\[-3mm]
  
Part (b) of the HvNW Theorem involves knowledge of 'almost crossings', i.e.,  of hyperbolic near approaches and almost touching  eigencurve pairs. These are currently identified  through  visual examination of the eigencurve plot such as depicted in Figure 2 for our example. We have not investigated how to accomplish this task computationally. Note that with $0 \leq t \leq 6$ in Figure 2,  the eigencurve pairs with labels 2 and 3,  5 and 6, 6 and 8, as well as 9 and 10  avoid crossings in a hyperbolic fashion and therefore each of these paired eigencurves must come from the same block in a diagonal block reduction of $A(t)$. If we enter their indices in  row-wise increasing ordered pairs in  the matrix 
$$Touch = \bp 2&3\\5&6\\6&8\\9&10 \ep  $$ 
 and invoke the command {\tt ve = almostTouch(ve,Touch)}, this will  adjust the previously computed  entries of $ve$ accordingly. For this example it generates
$ ve =  (1,    -1,    -1,     1,    -1,    -1,     1 ,   -1,    2,     2,    -2)$. Now we have refined the possible block decomposition for the matrix flow $A(t)$ from five blocks associated with  five labels 1, --1, 2, 3, and --3 to just four, namely  1, --1, 2, and --2. We might still be able to  improve our knowledge of the coarsest block-diagonal structure of $A(t)$ under similarities by looking at larger time intervals and thereby learn more about  additional  eigencurve crossings and/or  new almost  touching eigencurve pairs. 
\vspace*{-0mm}
\begin{center}
\includegraphics[width=120mm,angle=0]{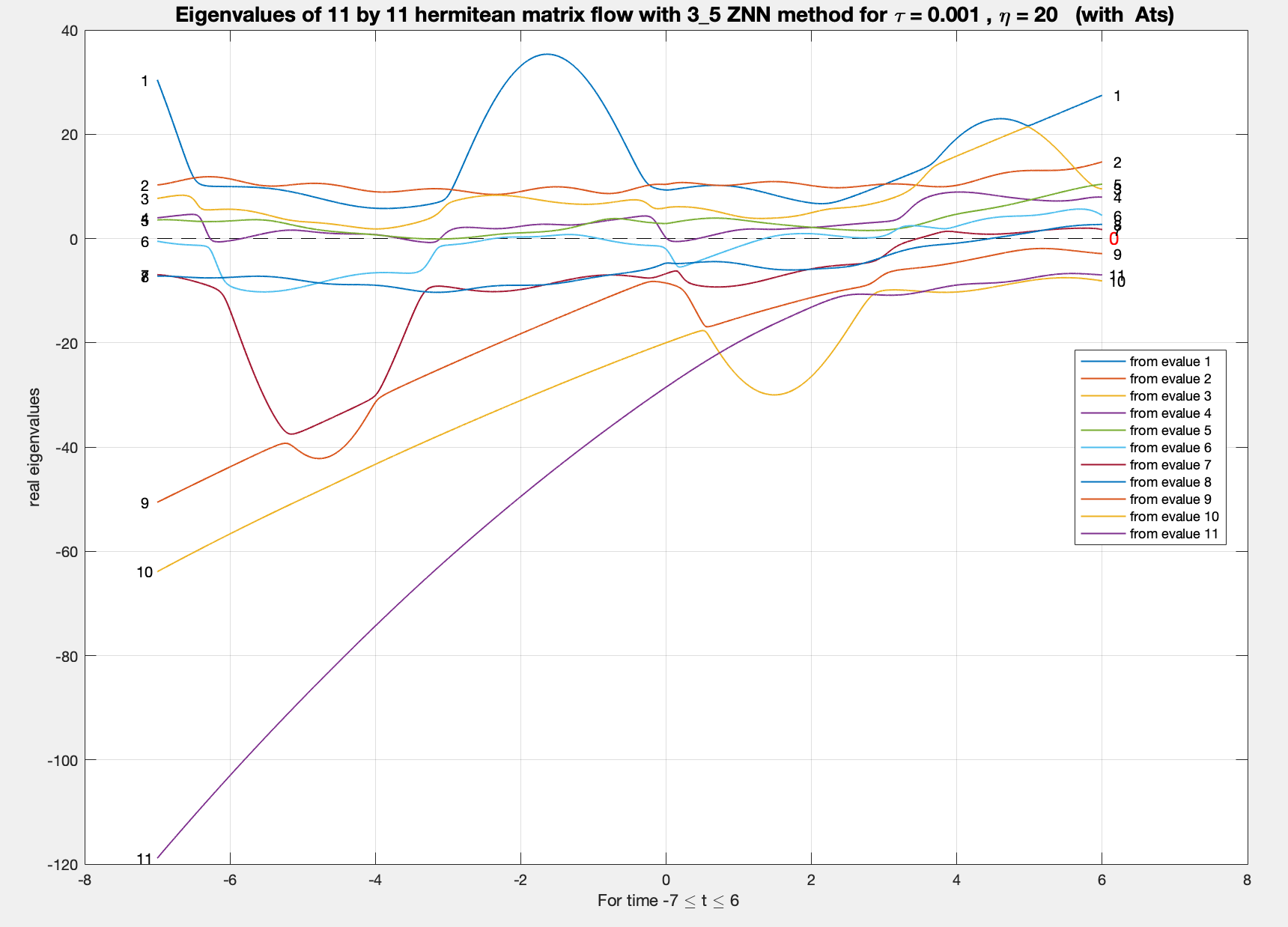} \\[0mm]
Figure  3
\end{center}
 The eigencurve plots  in Figure 3 for the same matrix flow $A(t)$ over the enlarged time interval $-7 \leq t \leq 6$ create  eigencurves  of $A(t)$ with six hyperbolic avoidances for the curves labeled 1 and 3, 3 and 4, 4 and 6, 6 and 7, 7 and 9, and 9 and 10. These label pairs are again stored in an 'almost touching' matrix $Touch$ for the enlarged interval. Calling  {\tt Chermitmatrixfloweig.m} for $A(t)$ and the revised time interval creates a new  eigencrossing data matrix R1 since the left time edge curve order obviously changes between $t_o = 0$ and $t_o=-7$ as new crossings occur on the leftwards extended time interval.\\[1mm]
\hspace*{56mm}\begin{tabular}{c|cccc}
curve & \multicolumn{4}{c}{crosses eigencurve}\\
number& \multicolumn{4}{c}{ with label}\\ \hline
1&2&0&0& 0\\
2&3&0&0&0\\
3&5&0&0&0\\
4&5&0&0&0\\
5&6&0&0&0\\
6&8&0&0&0\\
7&8&0&0&0\\
8&0&0&0&0\\
9&0&0&0&0\\
10&11&0&0&0
\end{tabular} \\[2mm]
From our new R1 and the new 'almost touching' matrix $Touch$,   the current coarsest decomposition into undecomposable diagonal blocks of $A(t)$ is indicated by the vector $ve = (1,    -1,     1,     1,    -1,     1,     1,    -1,     1,     1,    -1)$. Thus   $A(t)$ can be decomposed uniformly under one unitary similarity to  two irreducible diagonal blocks of dimensions 4 and 7, as the four +1  entries  and the seven --1  entries in $ve$ indicate.\\
Regarding speed and accuracy, when we run {\tt Chermitmatrixfloweig(11,0,6,3,5,0.001,20,1,1);} on the interval $[0,6]$, the ZNN computed eigenvalues agree in five leading digits with the eigenvalues of $A(t)$ at the right time interval end $t = t_f = 6$ when  we use the relatively low truncation error order convergent look-ahead finite difference formula {\tt  3\_5}. The computations take around 3 seconds. For the chosen sampling gap of $\tau = 1/1000$  our ZNN method   computes the complete 6,000 intermediate eigendata between $t_o$ and $t_f$ using essentially only one linear equations solver at each time step. For an explanation of ZNN's speed, its look-ahead difference formulas etc, see \cite{FUZhFoV,FU18}. How long an ODE path continuation method would take and how accurate it  would be when reaching $t_f = 6$, we do not know. \\[-2mm]

Can our dense complex hermitean matrix flow $A(t)$ be decomposed uniformly  into a direct sum of a 4-dimensional and a 7-dimensional undecomposable diagonal block? How was the matrix flow $A(t)$ constructed?  How can we construct a  coarsest block diagonalisation from a given dense hermitean matrix flow $A(t)$ somehow, if at all?\\[1mm]
 To form the dense 11  by 11 complex hermitean test  matrix flow $A(t)$, we started from a 7 by 7 single parameter variable  indecomposable general  matrix flow $B_1(t)$ with complex function entries and appended it to become  11 by 11 block-diagonal via the Matlab command  
 {\tt B2(t) = blkdiag(B1(t),2*B1(t)(2:5,2:5))}. Then we transformed the resulting general matrix flow $B_2(t) \in \CC^{11,11}$ to become the complex hermitean matrix flow  $B(t) = B_2(t) + B_2(t)^*$. And finally we obscured the  original block-diagonal form, comprised of one indecomposable 7-dimensional diagonal block and one indecomposable 4-dimensional one by replacing $B(t)$ with $A(t) = U^* B(t) U$ for a fixed complex random entry 11 by 11 unitary matrix $U$. The  matrix flow $A(t)$ that we used in our code thus  was  both complex hermitean  and dense. The matrices  $A(t)$  give no sign of decomposability. Yet our method has computed the indecomposable coarsest block diagonal structure  for $A(t)$ correctly by relying on the crossing geometry of its eigencurves and our code.\\
 Our  Matlab code {\tt Chermitmatrixfloweig.m} together with {\tt  almostTouch.m} and a bit of eye-balling for 'almost touching' hyperbolic crossings retrieved the hidden block-diagonal structure of this dense complex hermitean matrix flow readily once the time interval was chosen large enough. If we had eye-balled wrongly, {\tt  almostTouch.m} is designed to recognize such errors and it indicates which row in $Touch$ is erroneous. In this case we recommend to remove the row of offending  curve labels from $Touch$ and retry.\\[2mm]
Unfortunately we have no idea how to retrieve  the individual  underlying indecomposable 7 and 4 dimensional blocks of $A(t)$ themselves and find, for example, the generating decomposed matrix flow $B(t)$ (or a unitarily similar one)  for the given dense hermitean flow $A(t)$. All we can find right now are $A(t)$'s coarsest possible block dimensions but do not know how to compute its  individual diagonal blocks.\\[-3mm]

Our second example comes from a discussion of a 6 by 6 real symmetric test matrix flow and eigencurve crossings that appeared in 2010 on Mathematica's stackexchange.com site \cite{Ma2010} on "Tracking Eigenvalues through a Crossing". The poster, Jack S., suggests the symmetric matrix flow 
$$
B(t) =\bp 21t + 1/2     &       0&            0&           0&          0&            0\\
           0&   7t + 1/2&            0&           0& 7\cdot2^{1/2}t&            0\\
           0&            0&   1/2 - 7t&           0&            0& 7\cdot2^{1/2}t\\
          0&            0&            0& 1/2 - 21t&            0&            0\\
           0& 7\cdot2^{1/2}t&            0&           0&    14t - 1&            0\\
          0&            0& 7\cdot2^{1/2}t&           0&            0&  - 14t - 1 \ep
$$
and notices that 'Sometimes the eigenvalues may cross each other, but I want to make sure the right eigenvalue stays associated with its own state'. The site studies eigenvalue analyses and their potential tracking failures in parameter-varying matrix problems. ODE path continuing methods are offered, as well as characteristic polynomial approaches and a simple block diagonalization of the original matrix flow $B(t)$ followed by computing the parameter-varying eigenvalues for each separate diagonal block. Here we use ZNN and our algorithm again, but  on the camouflaged  real symmetric dense matrix flow $A(t) = U^TB(t)U$  for a random orthogonal entry matrix $U$.\\[1mm] 
 A call of {\tt Chermitmatrixfloweig(6,-.3,.1,5,6,0.0001,50,1,1)} takes 1 second and achieves 12 accurate leading digits for all computed eigenvalues of $A(t)$ at $t = t_f = 0.1$. The  eigencurve  crossings output is stored in  R below.\\[1mm]
\hspace*{56mm}\begin{tabular}{c|ccccc}
curve & \multicolumn{5}{c}{crosses eigencurve}\\
number& \multicolumn{5}{c}{ with label}\\ \hline
1&2&3&5&6&0\\
2&3&5&0&0&0\\
3&5&0&0&0&0\\
4&5&6&0&0&0\\
5&0&0&0&0&0\\
\end{tabular} \\[2mm]
The eigencurves are graphed in Figure 4.
\vspace*{-2mm} 
\begin{center}
\includegraphics[width=125mm,angle=0]{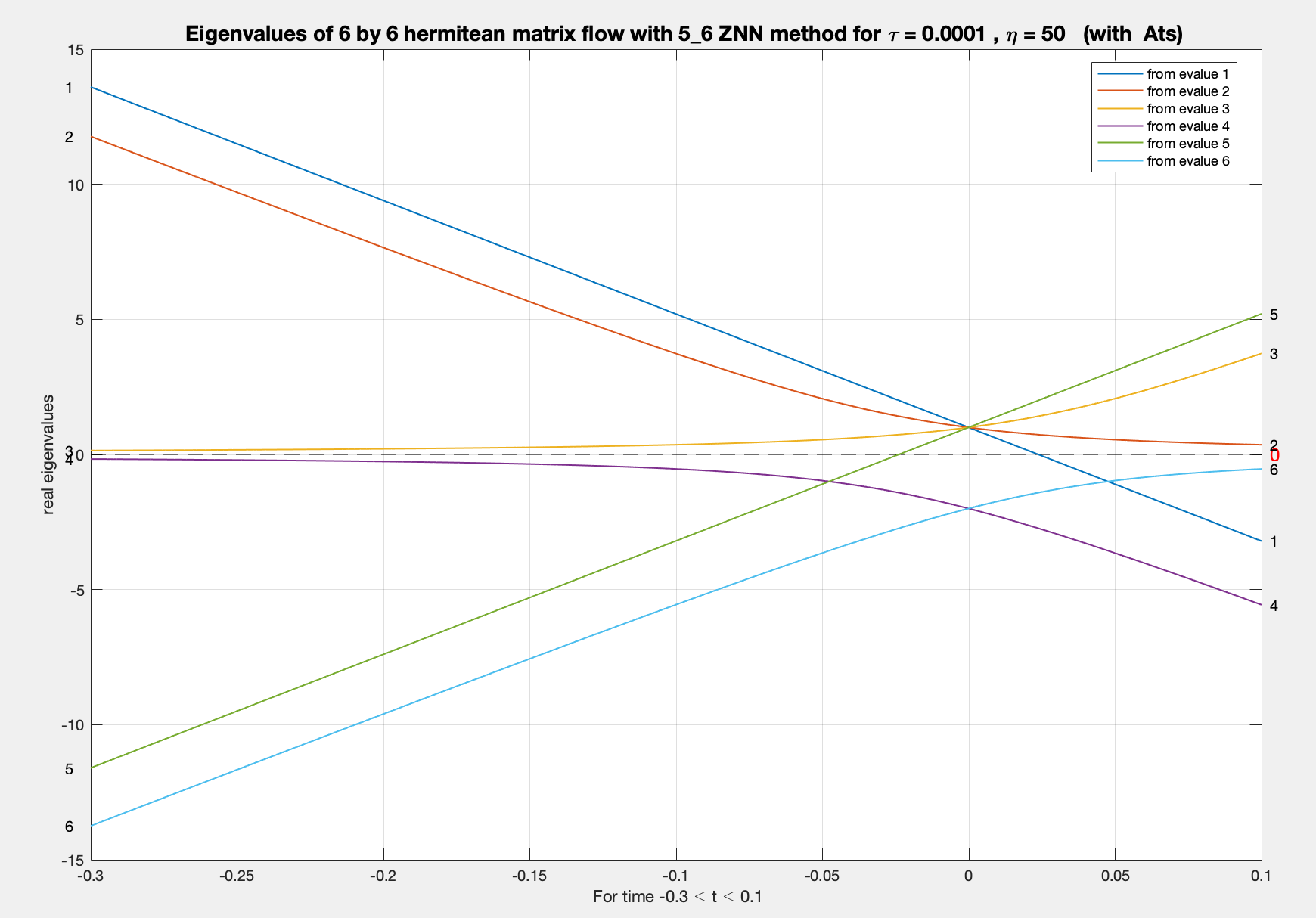} \\[0mm]
Figure  4
\end{center}
\vspace*{-1mm}
The above call computes $ve = (1,    -1,     2,     2,    -2,    -2)$. Note that the graph in Figure 4 shows no 'almost touching' eigencurve behavior. By looking at the label multiplicities in $ve$ we note that $A(t)$ can be reduced by an orthogonal similarity to a block diagonal matrix  composed of  four indecomposable  blocks, namely two of size 1 by 1  corresponding to the labels 1 and --1, and two  of size 2 by 2 for the repeated labels 2, 2 and --2, --2. This conforms well with $A(t)$'s origin in $B(t)$ which has precisely  this block-diagonal structure after applying a simple permutation similarity.\\[-3mm]

Finally we check our algorithm   with the  real diagonal 5 by 5 matrix flow seed \vspace*{-1mm}
$$B(t) =  \text{diag}(\sin(1-1/2t),1/2\cos(1/3t),\sin(t)\cos(-1-0.2t),\cos(2t-1/2),\cos(1+3t)^2)$$ 

\vspace*{-1mm}

\noindent
after it has been transformed by  an orthogonal random entry matrix $U_{5,5}$ into the dense real symmetric matrix flow $A(t) = U^TB(t)U$. Figure 5 shows no 'almost touching' hyperbolic evasions and dozens of eigencurve crossings.
 \vspace*{-1mm} 
\begin{center}
\includegraphics[width=110mm,angle=0]{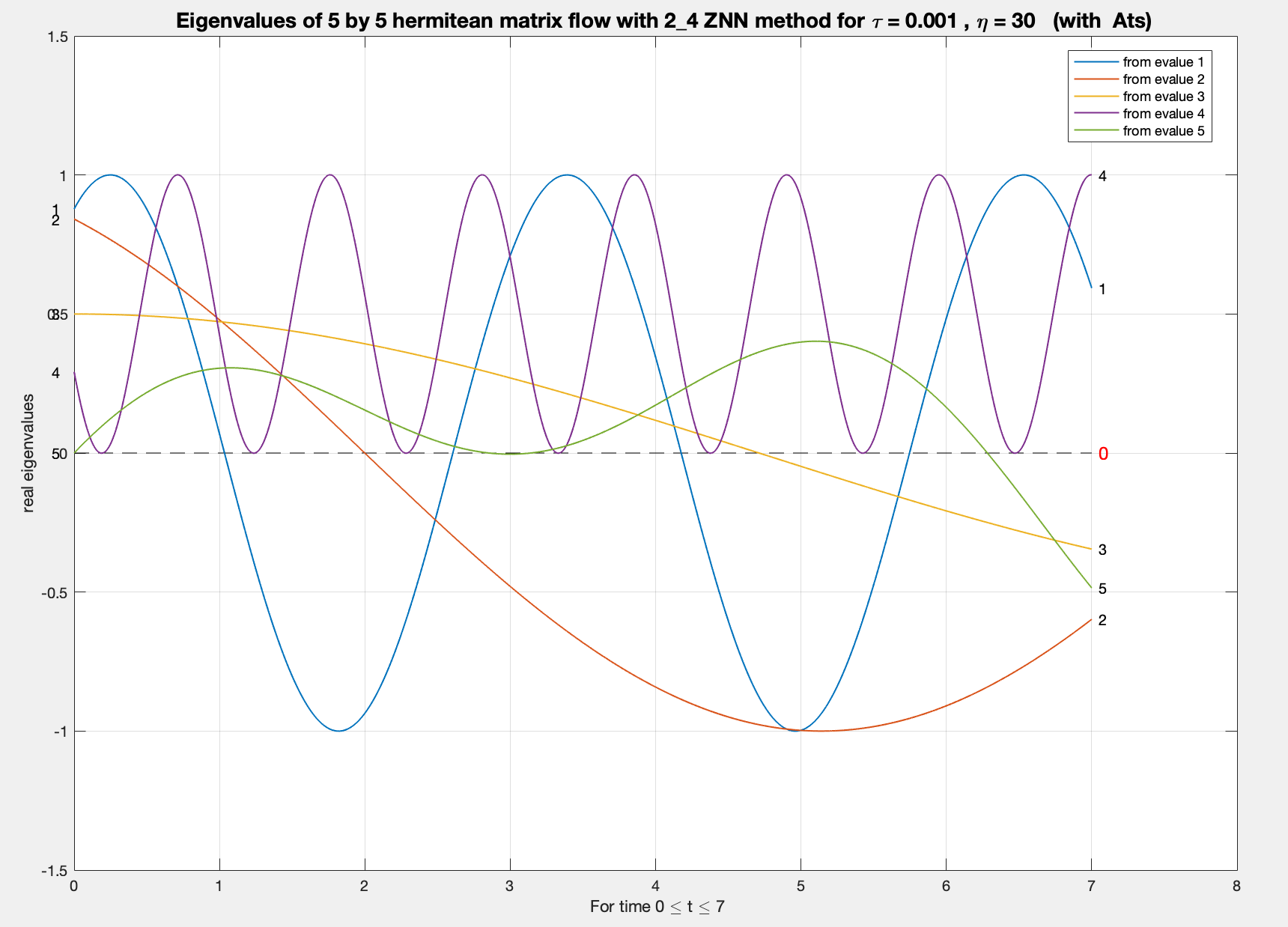} \\[0mm]
Figure  5
\end{center}
\noindent
The computed eigencurve crossing data matrix R1 and the block diagonalisation vector $ve$ for this example are\\ [1mm]
\hspace*{32mm}\begin{tabular}{c|ccccc}
curve & \multicolumn{5}{c}{crosses eigencurve}\\
number& \multicolumn{5}{c}{ with label}\\ \hline
1&2&3&4&5\\
2&3&4&5&0\\
3&4&5&0&0\\
4&5&0&0&0\\
\end{tabular}  \qquad and \qquad $ve =  \bp 1\\-1\\ 2\\-2\\3 \ep$ .\\[1mm]
Clearly {\tt Chermitmatrixfloweig} resolves this eigencurve and flow diagonalability problem perfectly in $ve$ which contains  5 distinct labels without repetition for our 5 by 5 symmetric flow problem, i.e., diagonalisation is the coarsest proper indecomposable block decomposition for $A(t)_{5,5}$.\\[-3mm]

We close the hermitean matrix flow eigencrossings section  with several open questions.\\[-3mm]

{\bf Question 2.1 :} What is the 'coarsest block diagonal form' of a hermitean matrix flow $A(t)$ depending on its eigencurve crossing geometry? Can our method compute the coarsest block diagonal form  reliably?\\[1mm]
As we realize, the given eigencurve crossings data may in fact be generated from an underlying  diagonal matrix flow $A(t)$ where every diagonal entry contains one respective eigencurve generating  function.  We would call such a diagonal decomposition the 'finest block-diagonalisation'. This may or may not be achievable from the given dense hermitean matrix flow $A(t)$. Instead, our algorithm tries to reduce the number of possible diagonal blocks in its computations by running through all  eigencurve crossings in turn. A large number of distinct entries in $ve$  signifies a relatively 'fine coarsest decomposition'. When the algorithm terminates after a potentially complete \emph{Touch} matrix incorporation, we believe or at least hope that we have gained insight into the 'coarsest block diagonalisation' of $A(t)$.\\[1mm]
 Are we really done then?  Can we be sure that this algorithm is complete in the sense that it can deal correctly with all possible eigencurve crossing data matrices  R1? Which matrices R1 are possible, which impossible to achieve as  eigencurve crossing matrices? How can these questions be answered mathematically and logically? Is knowing the eigencurve crossing matrix R1 sufficient to solve this problem? Of course not; the eigencurve crossing data matrix R1 is totally insufficient for indecomposable $n$ by $n$ hermitean matrix flows whose eigencurves are widely separated over the reals and never cross. This situation would make our algorithm's output vector $ve = (1,2,3,...,n)$ with $n$ distinct integers and indicate 1 by 1 block diagonalability in error. What other data would be needed in this worst case scenario? Do such matrix flows even exist?\\
  We do not know  answers to this set of questions.
  \newpage

{\bf Question 2.2 :}  Once we have found a proven correct decomposition algorithm for hermitean 1-parameter varying matrix flows $A(t)$ and know the coarsest block-diagonal structure that underlies $A(t)$, how can we find  a  conforming block-diagonal representation $B(t)$ of a given dense hermitean matrix flow $A(t)$ computationally? \\
 We have no answer to this question either.\\[-7mm]

\section{General 1-Parameter Complex and Real Matrix Flows and their Eigencurves}

\vspace*{-1mm}
In this section we depict and study general complex and real 1-parameter varying matrix flows $A(t)$ that are neither hermitean nor real symmetric. Such less restricted matrix flows generally give us 'wild and wooly' 3-D images in their eigencurve plots such as  Figure 6 below. This figure's  non-normal complex matrix flow  $C(t)_{11,11}$   was built from  our earlier complex 7 by 7 seed matrix $B_2(t)$ without the modifications to create a hermitean flow whose real eigencurves were   depicted in 2-D $\RR^2$ in Figures 2 and 3 earlier. Then $H(t) = U^*C(t)U$ is a dense general non-normal 11 by 11 complex matrix flow. To plot Figure 6 we use the general complex flow version  m-file {\tt Cmatrixfloweig.m} that is again based on the ZNN method for speed and accuracy. For general flows $H(t)$, $\RR^3$ describes the  parameter or time $t$ on one axis and the real and imaginary eigenvalue parts of $H(t)$ in the perpendicular plane at time $t$. When plotting time-varying complex eigencurves in $\RR^3$, it is very unlikely that 1-dimensional eigencurves will ever meet or cross -- unless, of course, the chosen complex matrix flow has a repeated block  such as $A(t) = \text{diag}(C(t),C(t))$ for example. 
\\[1mm]
For  'general' decomposable complex flows we  have never observed eigencurve crossings  which is standard with decomposable hermitean flows. From the depicted eigencurves of Figure 6, it appears impossible to assert whether the flow $A(t) $ allows a 4 by 4 and 7 by 7 block decomposition or not, this despite  its very creation as a decomposable flow. There is no crossing matrix to construct for generic non-hermitean complex matrix flows, or almost never,  since there literally appear to be no actual crossings in $\RR^3$ for complex 1-parameter general matrix flows. Our general complex matrix flow code {\tt Cmatrixfloweig.m}  checks on the minimal distances between individual eigencurves as a function of time $t$ for  decreasing distances of 1, 0.01, 0.0001, 0.000001 and shorter lengths. Even for examples with eigencurve distance minima  below $10^{-4}$ units of length we have never encountered a random entry complex flow example where two eigencurves got as close as $10^{-6}$ or $10^{-10}$ units.
\vspace*{-1mm} 
\begin{center}
\includegraphics[width=130mm]{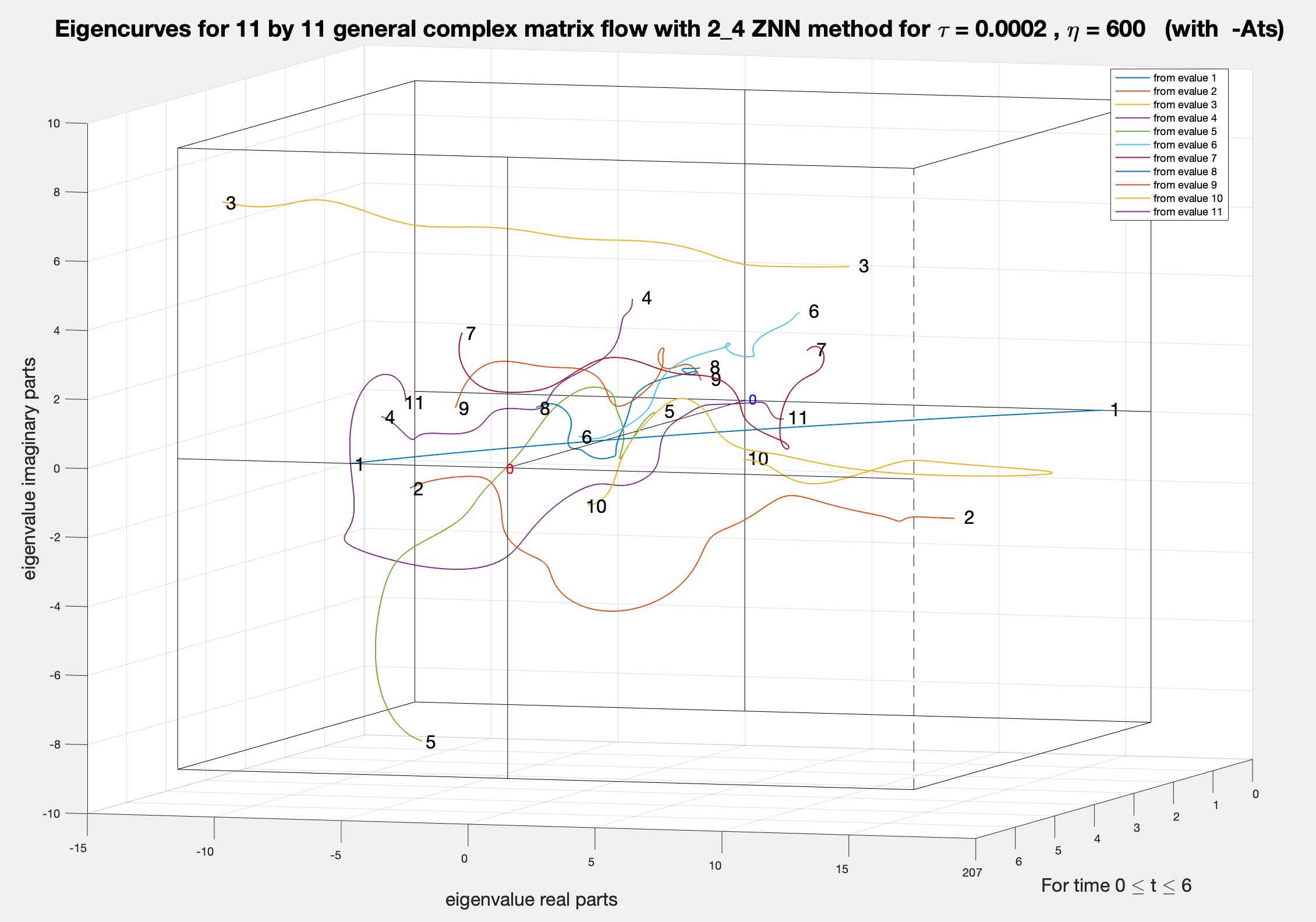} \\[0mm]
Figure  6
\end{center}
\newpage


How are the eigencurve geometry and the decomposability of general non-normal complex matrix flows related, if at all? We could find no literature and at first  had difficulties to  conceive of, let alone construct  suitable complex 1-parameter matrix flows whose eigencurves in $\RR^3$ might exhibit  the two crossing conditions (a) and (b) of the Hund-von Neumann-Wigner Theorem that is only known to hold for hermitean matrix flows.\\[1mm]
Eventually we were able to construct  a non-normal 1-parameter decomposable general 10 by 10 complex matrix flow $A(t)$ whose eigencurves exhibit HvNW like behavior. The set of Figures  7, 8, and 9  and further discussions below suggest that the two results (a) and (b) of Hund \cite{FH1927} and of von Neumann and Wigner \cite{NW} for hermitean matrix flows may hold for all general  complex and real 1-parameter matrix flows. A proof thereof is beckoning.\\[1mm]
This non-normal complex matrix flow  is composed of two indecomposable tridiagonal complex matrix flows and their block diagonal join. One of these is the 4 by 4 non-normal complex matrix flow\\[-8mm]

{\small 
$$
\hspace*{-10mm}A4(t) = \bp i(2-e^{t-1})+t/6 & 1                 & 0            & 0\\
                          1            & -2-2i\sin(t-1) & 1            &0\\
                          0            & 1                 &2i-2t           & 1 \\
                         0             & 0                     &1 & \sin(t+2)+it \ep ,  \ \ \ \text{\normalsize the second one   is 6 by 6 } 
 $$
 $$
A6(t) = \bp i-2\cos(2t)& 1                 & 0               & 0            & 0                                & 0\\
                          1            & -2-2i \sin(t-1) & 1            &0             & 0                                & 0\\
                          0            & 1                 &2i-t            & 1            & 0                               & 0\\
                         0             & 0                   & 1    & i e^{\sin(t)}       & 1                              & 0 \\
                         0             &  0                  & 0             & 1       &t/2+\sin(t)\cos(2it)/100 & 1\\
                         0             &  0                  & 0             & 0       &    1             & t-i/8\cos(it/3-1)\ep  
$$     
 {\normalsize and $A10(t)$ is the  concatenated  array } 
$A10(t) = \bp  A6(t)_{6,6} & O_{6,4}\\
                       O_{4,6} & A4(t)_{4,4} \ep \ .
$\\[-1mm] }

\noindent
Each of the above complex tridiagonal matrix flows is again made into a dense flow  via  one fixed unitary similarity transformation $B_{..} =  U_{..}^* A..(t)U_{..}$  before analyzing and plotting. Here the matrices $U_{..}$  have compatible dimensions 4 by 4, 6 by 6, or 10 by 10, respectively.\\
Figure 7 shows the eigencurves of $B4(t)$ when projected onto the eigenvalue real parts and time plane.
\vspace*{-2mm} 
\begin{center}
\includegraphics[width=120mm,angle=0]{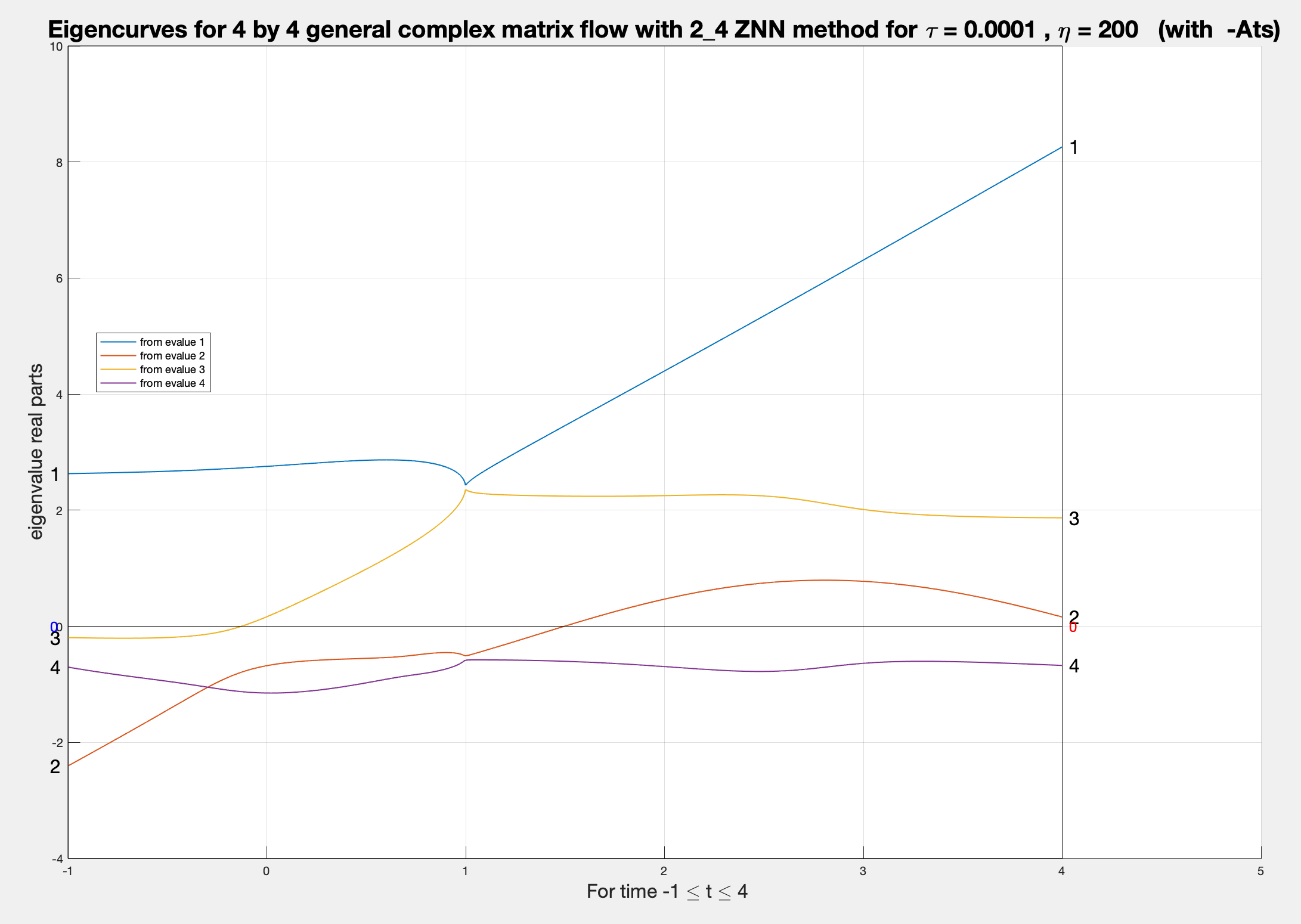} \\[0mm]
Figure  7
\end{center}
Figure 7  shows  'eigencurve avoidance' for the curve pairs 1 and 3 and 2 and 4 for the complex non-normal flow $B4(t)$ near $t = 1$.\\[1mm]
 Figure 8 below  shows a similar plot for $B6(t)$, again   with 'hyperbolic avoidance'  now near $t = 3.5$.
\vspace*{-2mm} 
\begin{center}
\includegraphics[width=120mm,angle=0]{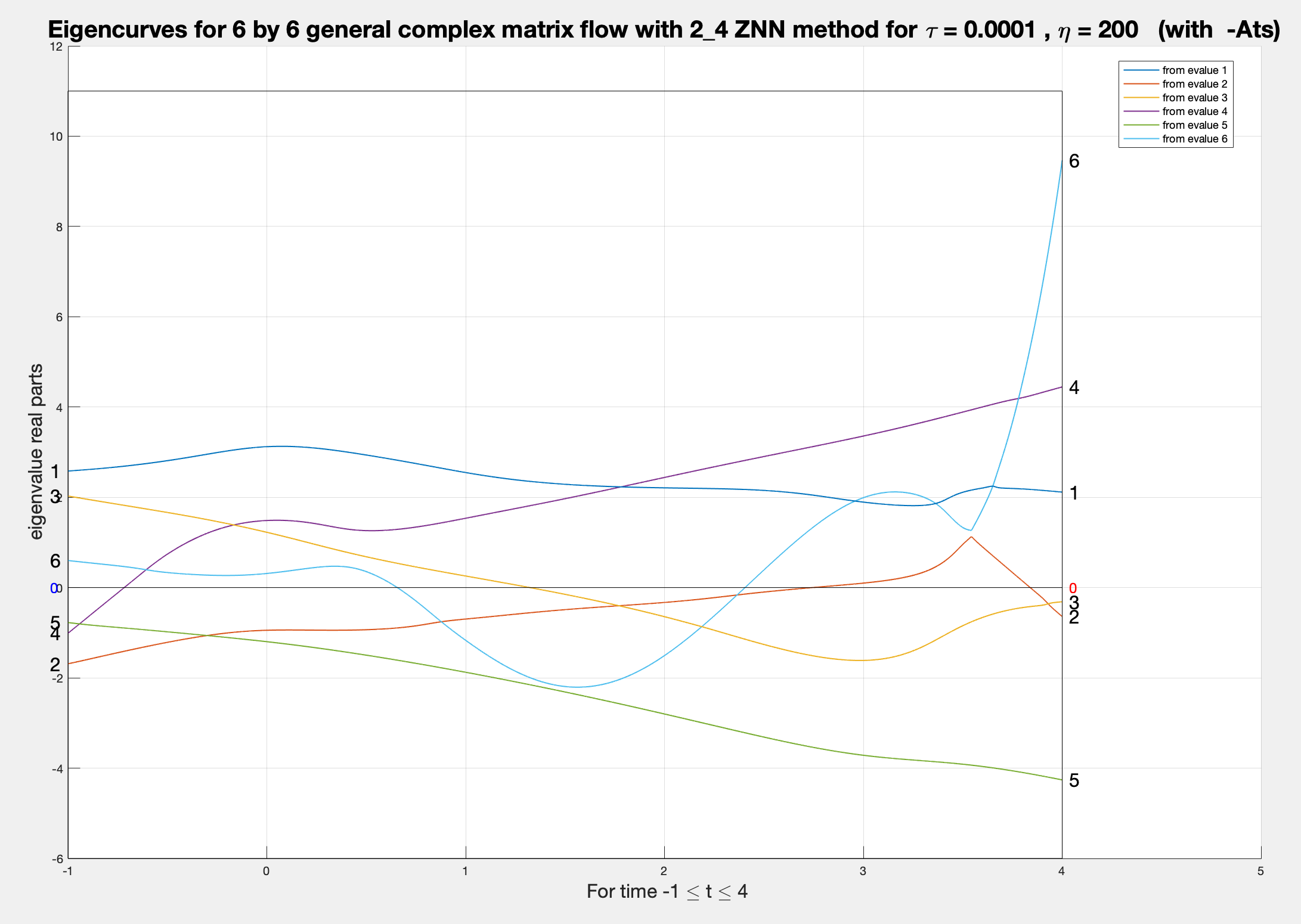} \\[0mm]
Figure  8
\end{center}
The combined, superimposed eigencurve  plot for the block-diagonal matrix flow $B10(t)$, again projected onto the real parts and time plane, is given in Figure 9 below.

\begin{center}
\includegraphics[width=120mm,angle=0]{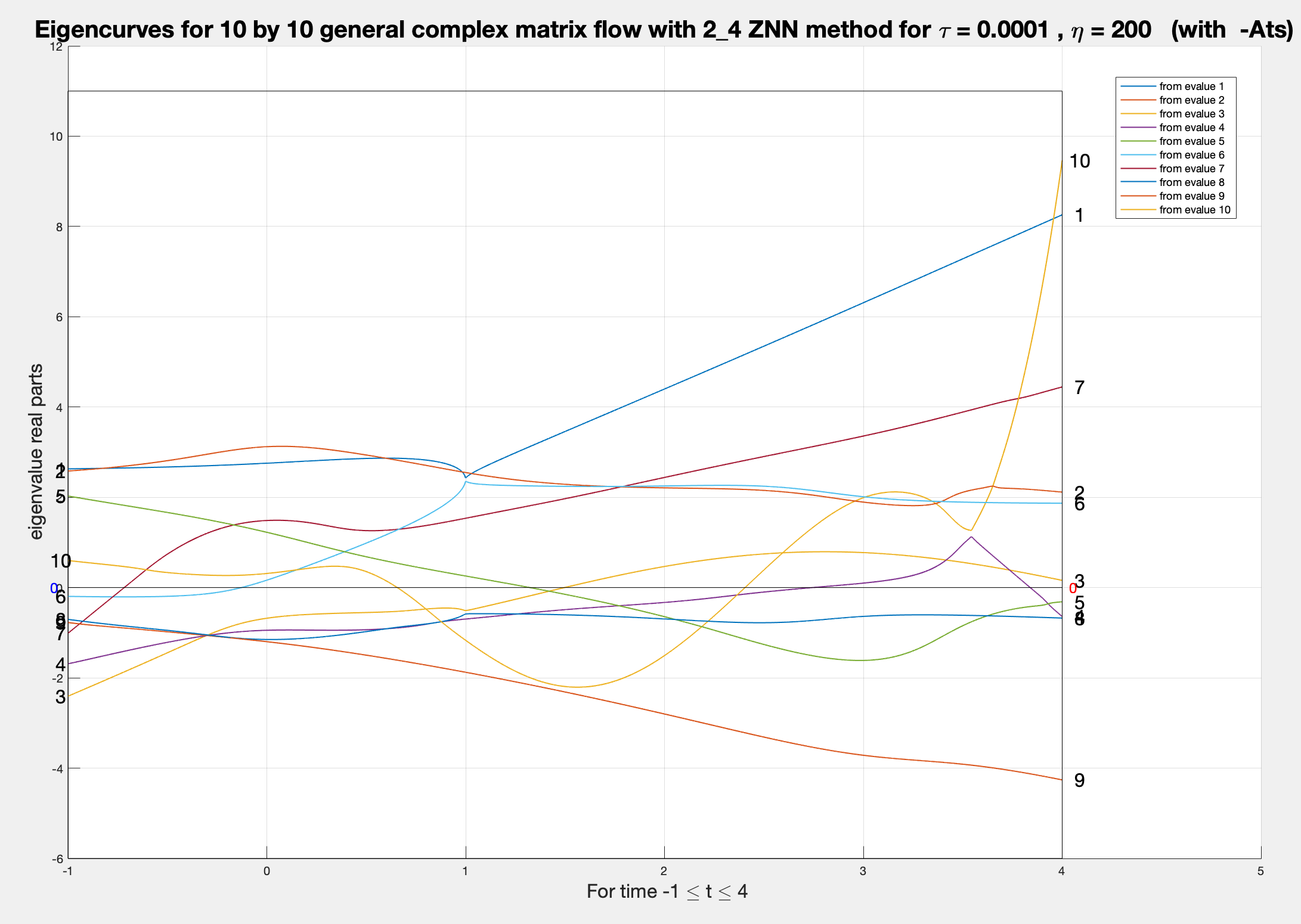} \\[0mm]
Figure  9
\end{center}
\vspace*{-2mm}
The eigencurve 'near touchings' for both $B4$ and $B6$ are  visible in the eigencurve plot for their concatenated flow $B10(t)$. The Rc eigencurve nearness data of $A10$ was computed after generating $B10(t)$ in {\tt formA10tricompl.m} via {\tt Cmatrixfloweig(10,-1,4,2,4,0.0001,200,1,-1)}. The nearness data of general complex flows is stored in Rc. The data in  Rc for our example indicates that the curves number 8 and 9 in Figure 9  get to within almost $10^{-2}$ units of each other, indicating that these two curves are close to crossing. But they show no signs of 'hyperbolic avoidance'. 
For help with reading the three plots in Figures 7 through 9, here is an equivalence list for the respective plot number labels and the eigencurves' 'almost touching' hyperbolic avoidance behavior.\\[2mm]
\hspace*{16mm}\begin{tabular}{lc|ccc|cl}
& \multicolumn{5}{c}{Equivalent Eigencurve Number Translations}&\\[1mm] 
 \hline \vspace*{-2mm}\\
& In Fig. 8 & \multicolumn{3}{c|}{In Fig. 9}& In Fig.7&\\
& for $B6(t)$& \multicolumn{3}{c|}{for $B10(t)$}& for $B4(t)$& \\[1mm]
\hline \vspace*{-2mm}\\
&&&1& $\rightarrow 4$&$1 \leftarrow 4$& alm. touch 6 \ --$|$\\
& $1 \rightarrow 2$& $1\leftarrow$& 2 &&& \hspace*{20.3mm} $|$\\
&&& 3&$\rightarrow 3$&$3 \leftarrow 3$&alm. touch 8 \hspace*{3mm} ---$|$\\
$|$-- \ alm. touch 6& $2\rightarrow 4$&$2\leftarrow$& 4 &&& \hspace*{20.3mm} $|$\hspace*{3.7mm}$|$\\
$|$ \hfill &$3 \rightarrow 5$&$3\leftarrow$&5&&& \hspace*{20.3mm} $|$\hspace*{3.7mm}$|$\\
$|$ \hfill &&& 6& $\rightarrow 6$& $6 \leftarrow 2$& alm. touch 1 \ --$|$\hspace*{3.7mm}$|$\\
$|$ \hfill &$4 \rightarrow 7$ & $4 \leftarrow$&7&&&\hspace*{25.9mm}$|$\\
$|$ \hfill &&&8&$\rightarrow 1$&$8 \leftarrow 1$&alm. touch 3 \hspace*{3.75mm}---$|$\\
$|$ \hfill &&&X&&&\\
$|$ \hfill & $5 \rightarrow 9$&$5 \leftarrow$&9&&&\\
$|$-- \ alm.touch 2& $6 \rightarrow 10$&$ 6 \leftarrow$& 10&&
\end{tabular}\\[2mm]
Near time $t_1 = 0.8174$ the eigencurves for $B10(t)$ with labels 8 and 9 pass each other in the real parts, imaginary parts, time space $\RR^3$ at $7.6 \cdot 10^{-3}$ distance units. Subtracting $ (.007.64534340607459 + .0001138407965141086 i) \cdot I_6$ from the first 6 by 6 tridiagonal block of $B10$, the eigenvalues of the resulting eigencurves with labels 8 and 9 of the  modified complex matrix flow $B10eps(t)$ agree in their 15 leading digits  at $t_1$ and for all practical purposes we have found a complex non-normal matrix flow with an eigencurve crossing. \\[1mm]
Looking at the imperceivably different eigencurve plot  for $B10eps(t)$  as displayed in Matlab's plot window, exhibits the three, now slightly shifted 'almost crossings'  of $B10$ clearly. And the two eigencurves labeled 8 and 9 now cross each other for $B10eps(t)$ and continue without interference, exemplifying HvHW's first condition (a) for general complex 1-parameter matrix flows and insinuating that the modified flow $B10eps(t)$ decomposes into at least two, possibly  indecomposable diagonal blocks. Note that  eigencurve 8 for the modified non-normal complex flow $B10eps(t)$ comes from $B4$ and eigencurve 9 derives from  the modified flow $B6$ according to our eigencurve number translations  list above.\\[1mm]
 Here, differing from the hermitean 1-parameter matrix flow case,  we have not been able to assert the actual number of coarsest blocks or their sizes.\\[-2mm]

An interesting  example from \cite{KST81} is the real non-normal  matrix flow
$$ A(t) = U^T \bp 1&t\\t^2&3 \ep U 
$$
for a fixed random entry real orthogonal 2 by 2 matrix $U$. $A(t)$'s two eigenvalues are complex conjugates for $t < -1$, double up as 2 and 2 at $t = -1$, and are distinct real for $ t > -1$, forcing its eigencurves to make right angle turns at $t = -1$.\\
 An accurate picture is plotted by the Matlab call of {\tt Cmatrixfloweig(2,-1.006,0,3,3,0.002,30,1,1)} in Figure 10 below for $-1.006 \leq t \leq 0$.
\begin{center}
\includegraphics[width=120mm,angle=0]{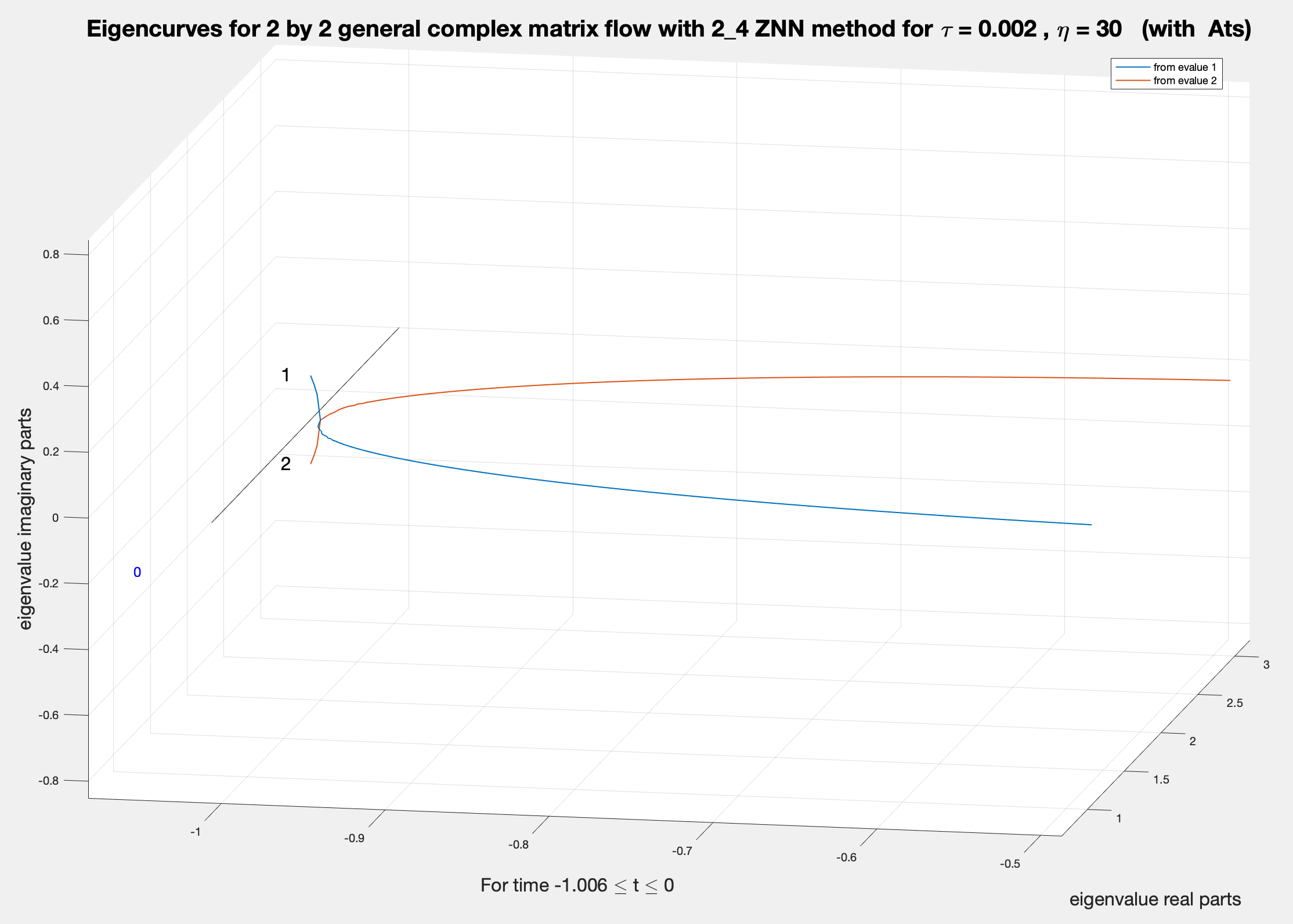} \\[0mm]
Figure  10
\end{center}
However, moving the starting time further back to $t_o = -2$ gives an incomplete picture.
 \begin{center}
\includegraphics[width=120mm,angle=0]{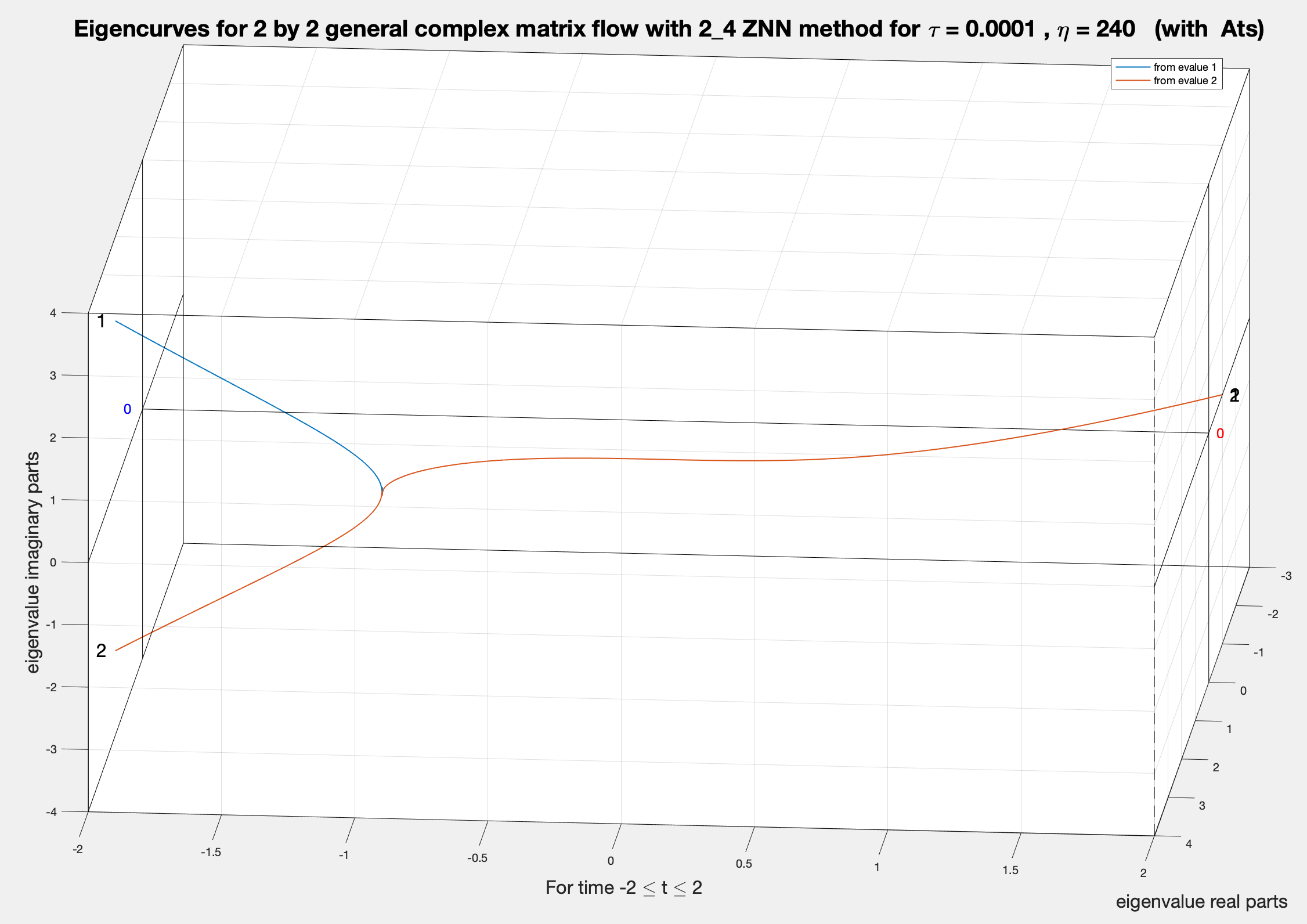} \\[0mm]
Figure  11
\end{center}
Figure 11 shows only one, now doubled eigencurve for $t > -1$ when computed via ZNN. We have observed similar behavior with eigencurve pairs for non-normal real matrix flows that transit from being complex conjugates to being real or vise versa. Sometimes a curve pair stays intact, coming in as a complex conjugate eigenvalue pair and leaving as a real pair. Sometimes the new eigencurve branches overlay and collapse into one and sometimes two formerly doubled-up complex conjugate or real eigencurve pairs proceed as two separate ones after such  90$^o$  degree direction flips. The reason for this behavior of ZNN eigen methods eludes us.   Similar glitches may occur with ODE IVP path continuation methods for non-normal real matrix flows or they may not. We are reminded here of Markus and Parilis's result in \cite{MP83} that  slight perturbations of a static matrix $A$ with a repeated eigenvalue $\lambda$ or repeated Jordan blocks for the same eigenvalue  $\lambda$ may only result in coarser, increased in size Jordan blocks for the perturbed eigenvalue $\tilde \lambda$ and never in  finer Jordan structures for the perturbed matrix $\tilde A$.\\[-2mm]

Based on the second 10 by 10 example of this Section, we conjecture a generalization of the Hund-von Neumann-Wigner Theorem \cite{FH1927,NW} :\\[1mm]
{\bf Generalized Hund-von Neumann-Wigner Theorem \ (a conjecture) \ :} (generalized HvNW) \\[1mm]
\hspace*{4mm}If $A(t)$ is an indecomposable general time-varying complex or real matrix flow, then\\[1mm] 
\hspace*{3mm} {\bf (a}) the eigencurves of $A(t)$ in the 3-dimensional real part, imaginary part and time space $\RR^3$ do not intersect.\\[1mm]
\hspace*{3mm} {\bf (b}) If two eigencurves of $A(t)$ approach each other, they veer off in a hyperbolic  way where the approaching\\
\hspace*{8mm}  space angle of either eigencurve  equals the leaving space angle of the other  after their close encounter.\\[-1mm]

Here are several related open questions.\\[-3mm]

{\bf Question 3.1 :}  Is the generalized HvNW Conjecture true? How can it be proved?\\[-3mm]

{\bf Question 3.2 :}  Given a 1-parameter varying complex non-normal and dense matrix flow $A(t)$, how can we assess its uniform decomposability from its eigencurve graphs in $\RR^3$? 
In 3-D space it is relatively rare for two  spatial curves to cross. This rarity  is due to the third degree of freedom  here when compared with  assessing the decomposability  of hermitean matrix flows from their eigencurve behavior which can be deduced entirely from $\RR^2$, see Section 2. How can we develop computational methods to find eigencurve crossings  and details of a flow's decomposability from its 3-D  eigencurve data? Is it enough, for example, to extend the parameter interval of the plots as done for hermitean flows?\\
We do not know.\\[-3mm]

{\bf Question 3.3 :} Are there any other criteria or data beyond eigencurve data  that can help  determine the uniform block-diagonalability for a general complex matrix flow and yield its block dimensions? What could they be?\\[-3mm]

{\bf Question 3.4 :} What lies behind the occasional doubling up of real matrix flow $A(t)$ eigencurve pairs that transit from being complex conjugate to real or in the reverse direction when evaluated via ZNN eigen methods?\\[-1mm]

{\bf Acknowledgement : } I am thankful to Nick Trefethen's mention of 'decomposing matrix flows' when I was visiting his group in Oxford in May 2019. His comment send me thinking differently and more deeply. It helped me to develop the  algorithms further and to formulate the conjecture.\\[-2mm]

{\bf Conflicts of Interest : }  There are none.

\vspace*{10mm}

\noindent
\centerline{{[} .. /box/local/latex/Coalescing Evalues.tex] \quad \today }

\bigskip

\noindent
11 image files :\\[2mm]
vNWhyperbola.png\\
herm11\_06.png\\
herm\_7to6graph.png\\
Math66sym.png\\[1mm]
5by5symdiag.png\\
complex11by11.png\\
Trinonnorm4by4.png\\
Trinonnorm6by6.png\\[1mm]
Trinonnorm10by10.png\\
2by2realdetail.png\\
2by2realwide.png


\medskip


\begin{thebibliography}{99}


\bibitem{DE99}{\sc Luca Dieci and Timo Eirola}, {\em
On smooth decompositions of matrices}, SIAM J. Matrix Anal. Appl., 20 (1999), p. 800-819 (1999).

\bibitem{DPP13} {\sc Luca Dieci,  Alessandra Papini and Alessandro Pugliese}, {\em Approximating coalescing points for eigenvalues of Hermitian matrices of three parameters}, SIAM J. Matrix Anal. Appl., 34 (2013),  p. 519-541. [MR 3054590], https://doi.org/10.1137/120898036 .


\bibitem{FH1927} {\sc Friedrich Hermann Hund}, {\em Zur Deutung der Molekelspektren. I.}, Zeitschrift f\"ur Physik, 40 (1927), p. 742 - 764.

\bibitem{KST81} {\sc Robert Kalaba, Karl Spingarn and Leigh Tesfatsion}, { \em Individual tracking of an eigenvalue and eigenvector of a parametrized matrix}, Nonlin. Anal., Theory, Meth. Appl., 5 (1981), p. 337 - 340.

\bibitem{LM} {\sc S\'ebastien  Loisel and Peter Maxwell}, {\em  Path-following method to determine the field of values of a  matrix at high accuracy}, SIAM J. Matrix Analysis, Appl., 39 (2018), p. 1726 - 1749.   https://doi.org/10.1137/17M1148608 .

\bibitem{MP83}{\sc A. S. Markus and E. E. Parilis}, {\em The change of the Jordan structure
of a matrix under small perturbations},  Lin. Alg. Appl., 54  (1983), p. 139 - 152.

\bibitem{Ma2010} 
{\url {https://mathematica.stackexchange.com/questions/165167/tracking-eigenvalues-through-a-crossing}}

\bibitem{NW} {\sc  John von Neumann and Eugene Paul Wigner}, {\em On the behavior of the eigenvalues of adiabatic processes}, Physikalische Zeitschrift, 30 (1929),  p. 467 - 470; reprinted in  {\em Quantum Chemistry, Classic Scientific Papers}, Hinne Hettema (editor), World Scientific (2000), p. 25 - 31.


\bibitem{Sch16}{\sc Kurt R. Schab, John M. Outwater, Matthew W. Young and Jennifer T. Bernhard}, {\em Eigenvalue crossing avoidance in characteristic modes},  IEEE Trans. Antennas Propag.,  64  (2016), p. 2617 - 2627, 
DOI: 10.1109/TAP.2016.2550098 .

 \bibitem{FU18} {\sc Frank Uhlig}, {\em The construction of high order convergent look-ahead finite difference formulas for Zhang Neural Networks},  J. Difference Equations and Appl., (2019), in print, 12 p.,  https://doi.org/10.1080/10236198.2019.1627343 .

\bibitem{FUZhFoV} {\sc Frank Uhlig}, {\em Zhang neural networks for fast and accurate computations of the field of values}, Lin. and Multilin. Alg., (2019), 14 p., in print, https://doi.org/10.1080/03081087.2019.1648375.

\bibitem{YZMeZNN} {\sc Frank Uhlig and Yunong Zhang}, {\em Time-varying matrix eigenanalyses via Zhang Neural Networks and finite difference equations}, Lin. Alg. Appl., 580  (2019), p. 417 - 435, https://doi.org/10.1016/j.laa.2019.06.028 .

\bibitem{FUEigencurves} {\sc Frank Uhlig}, The MATLAB codes for plotting and assessing  matrix flow eigencurves are  available at  {\url  {http://www.auburn.edu/~uhligfd/m_files/Eigencurves/}}

\bibitem{ZYLUH} {\sc Yunong Zhang, Min Yang, Chumin Li, Frank Uhlig, Haifeng Hu}, {\em New continuous ZD model for computation of time-varying eigenvalues and corresponding eigenvectors}, submitted, 16 p.

\end{thebibliography}
\end{document}